\documentclass{IEEEtran4PSCC}
\usepackage{xcolor}
\usepackage{booktabs}
\usepackage{amssymb}
\usepackage[normalem]{ulem}

\newcommand{\apriori}{HIS}
\newcommand{\nearrealtime}{NRT}
\newcommand{\notaddr}{\textbf{-}}

\ifCLASSINFOpdf
   \usepackage[pdftex]{graphicx}
\else
   \usepackage[dvips]{graphicx}
\fi
\usepackage{enumitem}
\usepackage{etoolbox}
\patchcmd{\thebibliography}{\section*{\refname}}{}{}{}
\usepackage{comment}
\usepackage{optidef}
\usepackage[]{amsmath}

\usepackage{tablefootnote}
\usepackage{algorithm}
\usepackage{algpseudocode}
\usepackage{dblfloatfix}
\usepackage{multirow}
\usepackage{hyperref}
\makeatletter
\let\old@ps@headings\ps@headings
\let\old@ps@IEEEtitlepagestyle\ps@IEEEtitlepagestyle
\def\psccfooter#1{%
    \def\ps@headings{%
        \old@ps@headings%
        \def\@oddfoot{\strut\hfill#1\hfill\strut}%
        \def\@evenfoot{\strut\hfill#1\hfill\strut}%
    }%
    \def\ps@IEEEtitlepagestyle{%
        \old@ps@IEEEtitlepagestyle%
        \def\@oddfoot{\strut\hfill#1\hfill\strut}%
        \def\@evenfoot{\strut\hfill#1\hfill\strut}%
    }%
    \ps@headings%
}
\makeatother

\psccfooter{%
        \parbox{\textwidth}{\hrulefill \\ \small{23rd Power Systems Computation Conference} \hfill \begin{minipage}{0.2\textwidth}\centering \vspace*{4pt} \includegraphics[scale=0.06]{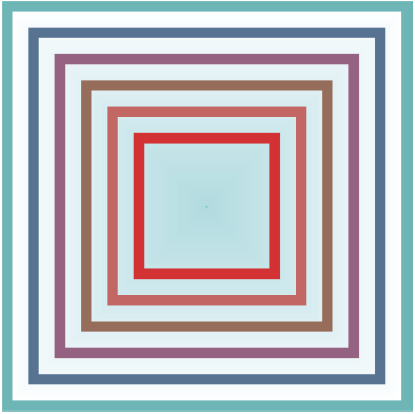}\\\small{PSCC 2024} \end{minipage} \hfill \small{Paris, France --- June 4 -- 7, 2024}}%
}
\begin{document}
\newcommand{\squeezeuptwo}{\vspace{-2mm}}
\newcommand{\squeezeupthree}{\vspace{-3mm}}

\title{Making  Distribution State Estimation Practical: Challenges and Opportunities}

%% To specify the authors when (number of affiliations <= 2)
% \author{
% \IEEEauthorblockN{Frederik Geth, Marta Vanin, Werner van Westering and Amritanshu Pandey}
% \IEEEauthorblockA{}
% }

\author{\IEEEauthorblockN{Frederik Geth \IEEEauthorrefmark{1},
Marta Vanin\IEEEauthorrefmark{2},
Werner Van Westering\IEEEauthorrefmark{3},
Terese Milford\IEEEauthorrefmark{4}, and
Amritanshu Pandey\IEEEauthorrefmark{5}
}
\IEEEauthorblockA{\IEEEauthorrefmark{1} GridQube, Brisbane, Australia,\\
\IEEEauthorrefmark{2} KU Leuven, Dept. Electrical Engineering (ESAT) and EnergyVille,
Leuven and Genk, Belgium,\\
\IEEEauthorrefmark{3} Alliander DNO, Arnhem, The Netherlands,\\
\IEEEauthorrefmark{4} Energex, EnergyQueensland, Brisbane, Australia,\\
\IEEEauthorrefmark{5} University of Vermont, Burlington, USA.
\vspace{-0.6cm}
}

}

%% To specify the authors when (number of affiliations > 2)
% \author{\IEEEauthorblockN{Author n.1\IEEEauthorrefmark{1},
% Author n.2\IEEEauthorrefmark{2},
% Author n.3\IEEEauthorrefmark{3}, 
% Author n.4\IEEEauthorrefmark{3} and
% Author n.5\IEEEauthorrefmark{4}}
% \IEEEauthorblockA{\IEEEauthorrefmark{1} Department Name of Organization A\\
% Name of the organization A,
% Address A\\ Emails if wanted}
% \IEEEauthorblockA{\IEEEauthorrefmark{2} Department Name of Organization B\\
% Name of the organization B,
% Address B\\ Emails if wanted}
% \IEEEauthorblockA{\IEEEauthorrefmark{3} Department Name of Organization C\\
% Name of the organization C,
% Address C\\ Emails if wanted}
% \IEEEauthorblockA{\IEEEauthorrefmark{4}Department Name of Organization D\\
% Name of the organization D,
% Address D\\ Emails if wanted}
% }

\maketitle

\begin{abstract}
In increasingly digitalized and metered distribution networks, state estimation is generally recognized as a key enabler of advanced network management functionalities. 
However, despite decades of research, the real-life adoption of state estimation in distribution systems remains sporadic. This systematization of knowledge paper discusses the cause for this while comparing industrial and academic experiences and reviewing well- and less-established research directions. 
We argue that to make distribution system state estimation more practical and applicable in the field, new perspectives are needed.
In particular, research should move away from conventional approaches and embrace generalized problem specifications and more comprehensive workflows. 
These, in turn, require algorithm advancements and more general mathematical formulations. 
We discuss lines of work to enable the delivery of tangible research.
%We review existing work that goes in this direction and highlight research 
 % \textcolor{blue}{plus small numerical examples?}.
\end{abstract}
\vspace{-0.5cm}
\begin{IEEEkeywords}
Distribution networks, Parameter estimation, State estimation, System identification
\end{IEEEkeywords}

% Use this to place sponsorships
% \thanksto{\noindent Submitted to the 23nd Power Systems Computation Conference (PSCC 2024).}

%\input{Introduction}

\section{Introduction}
% In the context of distribution networks (DNs), the energy transition entails impactful changes.
\noindent The distribution networks (DNs) in countries like Germany and The Netherlands are already feeling the impact of the energy transition, as they face congestion for hundreds to thousands of hours per year because of the rollout of rooftop solar, EVs, and loads like data centers \cite{Beckstedde2023}. 
This drives a need for real-time \emph{visibility} of the DN \emph{state} and congestions. 

\subsection{State Estimation and its Roots}\label{sec:state_estimation_roots}
\noindent In its standard textbook form~\cite{abur2004book}, power system state estimation (SE) is a statistical monitoring method that combines exact network models and noisy electrical measurements to calculate the \textit{most likely} state of the system, i.e., the minimum set of independent variables that fully describe the system's steady-state physics at a certain moment in time. 
% The conventional weighted least squares (WLS) form of SE is an optimal estimator under the (usual) assumption that all measurements are characterized by independent white (Gaussian) noise. 

%After unsuccessful attempts to compute bus voltages from manually collected measurements in the early 1960s~\cite{Wang2019}, 
The first real-time transmission system monitoring attempts directly read measurements from the supervisory control and data acquisition (SCADA) system~\cite{PajicPhD}. However, this had limited success due to individual measurements' noisy and unreliable nature. SE collates the different measurements, resulting in a reduced-noise, consistent and more reliable monitoring system that generates more information than the SCADA~\cite{PajicPhD}. SE established itself as the industrial practice in transmission system control rooms in the 1970s - following the cornerstone work by Schweppe et al.~\cite{Schweppe1}.%, Schweppe2, Schweppe3}. 

In its present industrial implementation, transmission system SE (TSSE) is close to the textbook form; it is static\footnote{i.e., time derivatives of state variables are not included in the model.}, and it presents a refresh rate of a few minutes. 
In addition to noise filtering, TSSE primarily 1) estimates variables that are not measured directly, 2) detects, identifies and handles bad data, 3) detects and identifies unreported topological changes (e.g., switch states), and 4) provides input states to several energy management system (EMS) functionalities, such as contingency analysis and real-time markets or economic dispatch.

While its adoption dates back to around 50 years ago, TSSE is still an active research field, and industrial implementations may  underperform. %, and, e.g., exacerbate extreme situations.  
% TSSE failure played a major role in the August 2003 Northeast blackout~\cite{PajicPhD}, and 
The NERC identified loss of TSSE %\footnote{mostly due to ``software" problems}
as the main source of EMS malfunction~\cite{NERC2020}. Gou and Shue~\cite{Gou2023} recently wrote on the challenges of industrial TSSE solutions. 

\vspace{-0.08cm}
\subsection{Changing Perspectives on SE in Distribution Networks} 

\noindent Distribution system SE (DSSE) literature only appeared about two decades after Schweppe's work~\cite{Baran1994}, when it became clear that distribution network (DN) management was becoming increasingly active and thus required monitoring. 

There has been general consensus that, given the differences between transmission and distribution systems, DSSE methods must differ from TSSE. 
Of the four primary use cases of TSSE defined in the previous section, {1)} still applies to DSSE, and is tied to the mathematical definition of SE. {2)} and {3)} are desirable, but existing TSSE techniques might not work for DSSE, given the scarcity of measurements and generally slow sampling rates. Moreover, {2)} and {3)} are very sensitive to errors in network models, which are common in DNs. 
These impact the SE residuals, making the handling of bad data and anomalous topology difficult if approached with TSSE-like residual analysis~\cite{abur2004book, Zanni2020}.   Use-case {4)} is of interest and possible for DNs. However, given the limited and `slow' measurements (which are also often averaged - as opposed to instantaneous), and few automated assets, the target EMS applications will differ from TSSE's, at least in the near future. 

\subsection{Data Quality Issues Demand Parameter Identification}\label{sec:dqi_demand_pi}
% Poor network data quality has been known to be problematic for TSSE for a long time. 
% For instance, Allemong reflects on a cultural challenge in organizations deploying TSSE \cite{Allemong2005}: ``It is probably a widely held belief that the equipment parameters, typically derived from system planning load flows, are unquestionably accurate, having stood the test of time.'' 

\noindent Network models for TSSE are considered trustworthy, %\footnote{Except for few (occasional) unreported topological changes.},
and  bad data detection and handling refers to \textit{measurement} data\footnote{Issues similar to bad data can be caused by cyber-attacks. While there is significant work on the topic, cyber-security is out of the scope of the paper.}. 
In DNs, the opposite is common: utilities tend to trust the (few) available sensors while network models are considered untrustworthy and error-prone. %\footnote{Trevizan et al.~\cite{Trevizan2021} show that indeed if multiple errors are present, meter readings are more accurate than voltage estimated by DSSE.}. 
The presence of errors in DN models is recognized both by industry and academia~\cite{GethCIRED2023, osti_1888157}. %

These network parameter errors affect DSSE in a fashion that is somewhat similar to bad measurements (i.e., they increase the residuals), but are more pervasive, persist over time (until the DN model is adapted) and  change over time (unreported system changes due to operations and maintenance). 
Conventional SE residual analysis struggles to differentiate between incorrect network model parameters and noisy measurements. A recent DSSE validation on a real DN illustrates that the large normalized residuals method may flag acceptable current measurements as bad data due to relatively modest impedance modeling errors~\cite{Zanni2020}.
As there are more potential error sources in DNs than the transmission system, calibrating DN models requires more expertise and time. 

Modern DN computation toolboxes, e.g., \textsc{OpenDSS}~\cite{OpenDSS}, \textsc{PowerModelsDistribution}~\cite{Fobes2020} and \textsc{PowerGridModel}~\cite{Xiang_PowerGridModel_power-grid-model}, support increasingly sophisticated models, e.g., explicit-neutral four-wire. 
However, their results are only as reliable as the input data used. 
The creation of DN models that accurately represent reality is needed to unlock the full potential of digital twins.

While conventional residual analysis might not always be effective, \textit{we posit that it is possible to exploit DSSE to build, improve and validate DN models}. Even more so, we believe that the maximal benefit from DSSE, from an industrial standpoint, rests in such model building tasks.
% \footnote{While other state-informed functionalities may still be possible and valuable as network automation increases, particularly in MV. }. %But DN models should anyhow be corrected beforehand.}. 
%, together with the joint estimation of states, parameter values, and unreported system changes\footnote{For instance switch changes which haven't been reported.}. 
This entails a paradigm shift, which consists of seeing state and parameter estimation as part of the same problem, as opposed to two disjoint problems, and is extensively discussed in this paper.

%In essence, this entails embedding state estimation and parameter identification in the same problem - as opposed to considering them two separate sequential problems. 
% This framework will allow us to construct up to date ``accurate" network models for distribution grids that have practical benefits for the DSO. 
% They can be used for hosting capacity, planning actions and other calculations, which are generally not accurate with model-free algorithms [REFERENCE]. 
% In particular, the latter appears as a non-solution given the DSOs' needs.

% Digital twins of power systems demand that the system being simulated matches the real-world one, which is a step beyond setting up a simulation model based on best practice.

\subsection{Scope and Contributions}
% To some extent, the canonical definition of what the knowns and unknowns are in (optimal) power flow ((O)PF)  is arbitrary; therefore, sticking to that definition in SE is also arbitrary. 
%States like transformer taps can be reasonably estimated based on a measurement snapshot with little risk of overfitting. Estimating states like line impedance, however, almost certainly leads to overfitting, unless  temporal trends are exploited. 

\noindent This work contextualizes and conceptualizes the problem space surrounding static DSSE, focusing on the applied industrial perspective. 
We observe that improvements needed in industry are far from mere pragmatic adjustments and require non-trivial theoretical and algorithmic advancements. 
We discuss the issue of making DSSE more practical by untangling some muddled terminologies, problem descriptions, and modeling approaches in the distribution system state and parameter estimation field.
%The ultimate goal is to inform researchers and practitioners of the gaps such that emerging  advancements are applicable in field. %Recommendations are also provided. 
We aim to inform and inspire researchers working on power system optimization methods to bridge the gaps and bring solutions to the field.

Alternative DSSE approaches exist, e.g., based on machine learning~\cite{Habib2023}, or dynamic filters~\cite{Uzunoglu2018}. 
Their discussion is out of scope as we focus on model-based static DSSE. %The use of synchrophasors is also not explicitly discussed. 

% These advancements maybe of interest to TSSE practitioners, too. Recently, the NERC \cite{NERC2020} identified loss of TSSE as the most significant source of EMS malfunction and `software' as the most significant driver. To this respect, Gou and Shue \cite{Gou2023} have recently written on the industrial challenges of established methods for state estimation in the transmission space, so this work complements that for the distribution context.

% We find the nomenclature in the existing system/parameter identification literature [REFERENCE] to be confusing and inconsistent, and as such dedicate Section~\ref{sec:terminology} to give an overview of existing terms and explain the terminological choices in this work. Among other, we want stress the differences between ``states" and ``parameters" and how different approaches are required to identify/estimate either.

The paper is structured as follows. 
First, we discuss the problem statements from the power flow (PF)/SE perspective compared to parameter estimation (PE) (\S \ref{sec:beyond_canonical_state_estimation}).
We review existing work on PE for DNs and find that it is often not integrated with DSSE (\S \ref{sec_state_param_estim}). 
% At the same time, we try to shed some light on the terminology used in the parameter identification field, which has developed in a  inconsistent and confusing manner. 
We then discuss the algorithmic advancements needed to enable practical and advanced DSSE (\S \ref{sec_algorithms}).
Next, we address practical challenges in DSSE deployment and the development of digital twins (\S \ref{sec_deployment}). 
Finally, \S \ref{sec_conclusions} summarizes and concludes the work.
%Finally, we conclude with unified vision for state and parameter estimation for distribution grids, highlighting future work directions (\S \ref{sec_conclusions}). 
% This is consistent with joint state-parameter estimation efforts on the transmission grids.

% \marta{other way to put the question: ``Making sense of bad distribution system state estimation results. Are the parameters or the data wrong?"}

\section{DSSE, related problems and misconceptions}
%\section{Misconceptions About State Estimation}
\label{sec:beyond_canonical_state_estimation}
\noindent We first discuss widely-accepted problem specifications of PF,  DSSE and related problems, and their inverses. 
% Next, we discuss the generalization to allow both to be performed simultaneously.

\subsection{Power Flow, State Estimation \& Load Allocation}

\noindent \emph{Power flow} is a mathematical problem that solves circuit physics while satisfying load/generation set points in active/reactive power given a reference bus. 
PF contains non-convex equalities from AC network constraints, necessitating nonlinear root-finding methods like the Newton-Raphson or fix-point iteration algorithms~\cite{Pau2023}. %\footnote{e.g. backward-forward sweep for radial networks~\cite{Pau2023}, current injection method for meshed.}. 

Conceptually, \emph{SE generalizes PF}\footnote{Any SE engine can solve a power flow by using the constant-power load/generation set points as measurements. 
However, a PF engine cannot be set up to solve nontrivial SE cases. }. 
Voltage measurements from any bus can be used directly in SE, whereas a PF can only accommodate complex voltage at the reference bus and voltage magnitude at generator buses. 
Other types of measurements can also be leveraged (e.g., relative and/or absolute voltage angles), as long as the corresponding functions define the measurement quantities in terms of the state variables. 
%Conventionally, the fitting criterion is a WLS of the measurement residuals.
% This results in a problem with non-convex quadratic constraints and quadratic convex objective, which can be solved with nonlinear WLS algorithms such as Gauss-Newton.

In the context of DNs, both PF and SE should model phase unbalance (and neutral voltage in four-wire DNs)~\cite{Melo2023,5937038}, and symmetrical components should be avoided \cite{5937038}. 
While the proposed modeling effort is not hard to support when building models and solvers from the ground up, re-purposing or adapting existing transmission-focused approaches and software solutions is practically intractable. 

\emph{Load allocation} (LA) is related to SE and PF, and typically conceived on top of a PF engine. 
LA aims to set up a realistic time simulation model of the system by assigning and dividing aggregated loads to a set of buses in the network~\cite{Kersting2008}. 
Practitioners then validate their model by looking at voltage and power at the substation but typically have no other measurement data available to assess the LA quality. 
% LA problems can be set up using SE engines, and provided that the system is observable without measuring the loads (e.g., thanks to voltage measurements), these can be best allocated by SE than other tools. %The key is the availability of sufficient voltage measurements.

% Note that practitioners invariably make assumptions regarding the voltage-dependency of loads. 
% The assumption of constant power behavior is very common, as it represents the worst case w.r.t. voltage drops. When more realistic assessments are required, loads are additionally parameterized (partially) by constant current, constant impedance or `exponential' load models~\cite{Claeys2021}.
% In the context of SE, these nuances are typically neglected but could be incorporated. For instance, including load information using power or ZIP models are likely to result in mathematically different outputs.

\subsection{States vs. Variables vs. Parameters}

\noindent Conventionally, for the (static) power system SE problem, complex nodal voltages are taken as the only unknowns and are represented in a stacked fashion in the `state variable  vector'. %, which is composed of complex voltage variables only. 
This restrictive definition excludes transformer taps (and other parameters). However, taps may change in the field based on local control loops, %(e.g., a hysteresis controller looking at voltage magnitude). 
and thus their estimation/tracking may improve the quality of the results. 
% This argument is easily extended to other parameters. 
Table \ref{tab:categories}  summarizes the features of relevant PE problems.

\begin{table*}[b]
\vspace*{-\baselineskip}
\caption{Parameter estimation categories.  }
\centering
\begin{tabular}{l l p{9.5cm} l } 
\hline
Category & Specific parameters & Problem statement & Nature   \\
\hline 
Topology  & Topology detection & Obtain list of energized lines/cables given a trustworhy graph & discrete  \\
 & Topology estimation (green-field) & Given a set of connected buses, estimate the graph & discrete \\
 & Topology estimation (incremental) & Given a graph, fix the connectivity errors. Includes switch state estimation. & discrete \\
 % & Switch state estimation & Same as topology detection & discrete \\
 & Phase connectivity identification & Estimate phase connectivity, usually given meter-to-transformer/feeder assignment & discrete \\
 & Meter-to-transformer assignment & Find to which transformer a user/meter is connected & discrete  \\
\hline
 Line  & Impedance estimation & Infer line impedance matrices & continuous * \\
 & Line type estimation & Select most likely construction code from a set of construction codes& discrete \\
 & Line properties estimation & Estimate most likely temperature,  line spacings, radius, material types & mix \\
 % & Carson's equations estimation & Estimate most likely line spacings, radius, material types, etc. & mix *\\
 \hline
 Transformer  & Impedance estimation & Estimate transformer loss terms & continuous \\ 
 & Transformer tap position & Estimate tap position of the transformer & discrete \\ 
% & Transformer phase shifter  & Estimate tap phase position & discrete &snapshot\\
 & Transformer type & Estimate the type of transformer (wye-wye, delta-wye, center-tap) & discrete \\
 \hline
 Regulator /   & Gang vs per-phase operation & Categorize whether taps are locked together or not & discrete \\
 OLTC & Impedance estimation& Estimate impedances in equivalent circuit& continuous  \\ 
   & Controller mode estimation & Categorize hysteresis vs line-drop compenstation, etc. & discrete \\
  & Bidirectional capable? & Estimate whether controller operates correctly under reverse flows & discrete  \\ \hline
  PV system   & Rating estimation & Estimate ratings of panels and inverter& continuous  \\
  & Control mode identification & Categorize constant PF vs volt-var/watt & discrete  \\
  & Control setting estimation & Estimate volt-var/watt curve values & mix \\
 \hline
\end{tabular}
\label{tab:categories}
\vspace*{-1.5\baselineskip}
\end{table*}

In typical DSSE implementations, the Jacobian (and Hessian) derivative expressions are hard-coded, and adding new state variables requires adaptation of these matrices, which is error-prone. 
Relying on modern optimization frameworks that support automatic differentiation enables faster prototyping. 
However, much of the current research in power grid state estimation does not consider progress being made in the applied mathematics or operations research community. 
In this context, we observe that a majority of the leading teams of the \emph{ARPA-E Grid Optimization Competition} resort to the same or similar underlying NLP algorithm \emph{implementations}~\cite{Safdarian2022}, e.g. the open-source solver \textsc{Ipopt} or KNITRO.

% \footnote{Recent optimization software frameworks do allow for \emph{automatic differentiation} though, enabling faster prototyping.}.  

%\cite{ZhangPESGM2015} recovers (Single-phase) transformer parameters in real time (for state-of-health): turns ratio, series winding resistance, series leakage inductance, shunt magnitude, shunt core loss resistance. 
%\cite{Aravind2021} estimates parameters for steady-state operation using transformer terminal measurements.

% \begin{table*}[b]
% \caption{Parameter estimation literature. N/A means that it doesn't make sense.}
% \centering
% \begin{tabular}{l l l l l} 
% \hline
% Parameters & Joint to state & Disjoint & Sequential & \apriori \ or \nearrealtime  \\
% \hline 
%  Topology detection &  \checkmark: & \\
% Topology estimation (green-field) &  & \\
%  Topology estimation (incremental) & & \\
%   Switch state estimation & \\
%  Phase connectivity identification &  & \\
%  Meter-to-transformer assignment & & \\
% \hline
%  Impedance estimation & &  \\
%  Line type estimation & & \\
%  Line temperature estimation & & 
%  \\
%  Carson's equations estimation & & \\
%  \hline
%  Impedance estimation & & \\ 
%  Transformer tap position &  & \\ 
%  Transformer phase shifter  & &  \\
%  Transformer type & & \\
%  \hline
%  Gang vs per-phase operation &  & \\
% Impedance estimation & & \\ 
%    Controller mode estimation & & \\
%   Bidirectional capable? & & \\ \hline
%   PV rating estimation & & \\
%   PV control mode identification & & 
%  \\
%   PV control setting estimation & & \\
%  \hline
% \end{tabular}
% \label{tab:parameter_literature}
% \end{table*}

\subsection{Forward and Inverse Problems}\label{sec:forward_inverse_problems}

\noindent In the literature, parameter estimation problems for power networks have been discussed under various labels: 
\begin{itemize}[leftmargin=*]
    % \item Parameter estimation
    \item system or parameter \emph{identification} or \emph{estimation};
    \item \emph{generalized} or \emph{augmented} state estimation;
    % \item state augmentation;
    \item \emph{inverse power flow} problem.
\end{itemize}
In the literature, the terms are used interchangeably,  e.g.~\cite{Kurup2023, Bariya, Wang2020}. Here we use `parameter estimation' as the overarching term to indicate the  concept. For individual parameters, we stick to the continuous/discrete considerations. 

% We choose the language of inverse problems to delineate our work. 
To minimize semantic confusion, we use terminology  inspired by the language in the operations research community:
\begin{itemize}[leftmargin=*]
    \item parameters: representing knowns  in the mathematical model, i.e.,  inputs supplied to the engine;
    \item variables: representing unknowns in the mathematical model, i.e.,  results calculated by the engine.
\end{itemize}
% Semantically, \textit{estimation} suggests an approximate assessment, and appears better suited for continuous parameters. \textit{Identification}, on the other hand, suggests a labelling effort, which appears more suited for discrete parameters. 

Mathematically speaking, the inverse problem of a (forward) problem consists of recasting all parameters as variables and vice-versa in the postulated system of equations. 
However, most authors only involve the continuous parameters/variables in this swapping operation, as making discrete parameters into variables may trigger the problem to become combinatorial in nature.
As PF and SE are generally considered continuous problems, when postulating their inverse, the discrete inputs - i.e. the topology - are typically kept as-is.
Nevertheless, this restriction can be dropped, with authors investigating discrete topics like topology identification (TI).

Table \ref{tab:stateestimationarchetype} summarizes the variables and parameters as defined in canonical PF and SE problems, and illustrates continuous and continuous+discrete generalizations. 
Potentially, \emph{all} symbols appearing in the equations of the mathematical model may become variables. In such case, the only model parameters (inputs) are the measurements, and we end up intersecting SE, PF and OPF with their own inverses.
%As we are generalizing the problems here, eventually all symbols appearing in the equations of the mathematical model may become variables, and we end up intersecting SE, PF and OPF with their own inverses.
However, the challenge of exploiting the mathematical structure where all symbols are variables is that the overfitting risk increases. Authors that address joint DSSE and PE often assume a fixed topology (except perhaps a few variables), and thus the problems retain largely continuous features, whereas TI problems abstract on continuous variables so that they are purely combinatorial.

\begin{table*}[b]
\vspace*{-\baselineskip}
\caption{Canonical problem specifications and beyond}
\centering
\begin{tabular}{l l l l l l l l l} 
\hline
&         &        & \multicolumn{2}{c}{Canonical power flow}                     &\multicolumn{2}{c}{Canonical state estimation} & Generalized cont. & Generalized\\
         \cmidrule(lr){4-5}  \cmidrule(lr){6-7}   \cmidrule(lr){8-8} \cmidrule(lr){9-9} 
&        Quantity & Symbol & Forward & Inverse & Forward & Inverse & Combined  & Combined \\ 
\hline
Discrete &        Line connectivity & $(l,i,j) \in \mathcal{T}$ & param. & param.& param. & param.  & param. & var.\\
&        Load connectivity & $(d,i) \in \mathcal{T}^{d}$ & param. & param.& param. & param.  & param. & var.\\
&        Generator connectivity & $(g,i) \in \mathcal{T}^{g}$ & param. & param.& param. & param.  & param. & var.\\
\hline 
Continuous &        Bus voltage & $U_{i}$ & var. & param.& var. & param. & var.& var.\\
&        Line impedance & $Z_l$ & param. & var.& param. & var.& var. & var.\\
&        Load value & $S_d$ & param. & var. &  var. & param. & var. & var.\\
&        Generator value  & $S_g$ & param. & var. &  var. & param. & var.& var.\\
% &       Measurements & & - & - & param. & param. & param.  & param.  \\
% \hline
% Derived &        Line current & $I_{lij}$ & var. & param. & var. & param.  & var. & var.\\
% &        Nodal adm. matrix & $\mathbf{Y}$ & param. & var.  & param. & var.& var. & var.\\
\hline
Extensions&        TF tap ratio & $T_l$ &param. &param.&param. &param. &  var. & var.\\
        & Those listed in Table \ref{tab:categories}& &- &- &- &- &  var. & var.\\
        & Variable bounds & & -& -& - & - & $\checkmark$& $\checkmark$ \\
 \hline
Math. & Class&&Nonl. root-find.& NLP &Nonl. WLS& NLP  &NLP & MINLP \\
complexity & Algorithm&&NR & IP &GN& IP &IP & custom \\
 \hline
\end{tabular}
\label{tab:stateestimationarchetype}
\vspace*{-1.5\baselineskip}
\end{table*}

We note that Yuan et al. \cite{9858017} defined the inverse PF problem as: ``the estimation of the nodal admittance matrix from synchronized measurements of voltage and current phasors''. As per our mapping, that would be more precisely described as the inverse \emph{state estimation} problem.  Note that the generalization \emph{hierarchy reverses}: SE generalizes PF in the forward space, and viceversa in the inverse space. %Note that for inverse problems, the generalization hierarchy is reversed: inverse PF generalizes inverse SE. %whereas (forward) SE generalizes (forward) PF. 
%We contrast this with Yuan et al. \cite{9858017} defining the inverse PF problem as: ``the estimation of the nodal admittance matrix from synchronized measurements of voltage and current phasors''.
%According to our mapping, that could be more precisely described as the inverse \emph{state estimation} problem.  Note that for inverse problems, the generalization hierarchy is reversed: inverse PF generalizes inverse SE whereas (forward) SE generalizes (forward) PF. 

The scope of any problem can be increased, e.g., by adding new variables (like taps in SE), or by allowing new mathematical structures such as variable bounds. 
Once an engine is built for a stated problem specification, it is often feasible to re-cast variables as parameters, while the opposite would increase the complexity and might be impractical.

%By intersecting the forward and reverse problems, in principle, all symbols appearing in the system of equations become variables, such that the only model parameters (inputs) are the measurements. 

%Note that many authors of joint state and parameter estimation still assume the topology (except for a few switch states, possibly) and phase connectivity fixed and as inputs, and therefore the problem largely retains a continuous structure\footnote{Except for properties such as transformer taps, which are discrete.}. 
%Once topological features are considered variable or unknown, the problem is necessarily combinatorial. 

% The \emph{nodal admittance matrix} offers a solid basis to look at uniqueness properties \cite{9858017} .
% However, while theorically useful, it is insufficient when looking at the problem of cleaning up existing network data sets \emph{as they are stored in the utility databases}, as nodal admittance matrices have to be rebuilt when the network changes.
% \amrit{Anyway, can we formally show the disadvantages of characterizing the network with just a nodal admittance matrix? Maybe, show the benefits of alternative approaches?} 

When PE is performed, it is often crucial to exploit time series, so the generalized problems need to be conceived of as \emph{fundamentally multiperiod}: assuming  parameter values to be time-independent helps  manage the overfitting risk.

\vspace{-0.1cm}
\subsection{Uniqueness and Observability}

\noindent From a  mathematical/control theory standpoint, observability for a steady-state system $S$ requires a one-to-one mapping between measurements and states in a noiseless environment. Given a measurement vector $z$, states $x$, and noise vector $n$ for some system $S$, observability holds if, for a noiseless scenario $n=0$,  we can recover $x$ from $z$ uniquely.
Local observability is less strict and can be met if, within a small operating neighborhood of system $S$, there is a unique   $x$ to $z$ mapping in a noiseless scenario. Full-rank conditions, as generally used in power systems SE, can  satisfy local observability.

Uniqueness, in the context of root-finding, requires a unique mapping between an arbitrary system of nonlinear equations $F(x)$ and $x$, in some interval $(\overline{x}, \underline{x})$.
Uniqueness issues exist:
\begin{itemize}[leftmargin=*]
    \item  PF generally does not have a unique solution, but in some practical intervals of states this can be the case
 (e.g., high-voltage solutions).
    \item SE is unique in observable linear systems. In nonlinear ones, generally only local observability can be achieved.
    \item Sometimes SE is only locally partially observable; here, only some states $x_o \subset x$ can be estimated given a set of limited measurements $z_o$.
    \item LA is generally not unique. If performed with a SE engine, the SE considerations hold.   
    \item Nodal admittance 
 matrices are usually unobservable; careful construction achieves local observability~\cite{9858017}. Other than~\cite{9858017}, observability and uniqueness are rarely formally discussed for PE.
\end{itemize}
The industrial view is that uniqueness is not that interesting, \emph{unless} it leads to wrong outcomes, i.e., losing extrapolation value. 
Out-of-sample testing is a good strategy to build confidence in models that are not guaranteed to be unique.

\subsection{Error and Uncertainty Sources}

%In the distribution network today, it is hard to distinguish various sources of erroneous data.
%Furthermore, the cost of investigation of issues may be prohibitive.
\noindent Existing implementations of DSSE and PF platforms suffer from  inaccuracies across different data sources~\cite{GethCIRED2023}:
\begin{itemize}[leftmargin=*]
   % \item \emph{modeling errors} or shortcuts due to applying the circuit laws under invalid assumptions, e.g., using Kron's reduction (KR) in networks where neutrals aren't pervasively grounded;
   \item \emph{modeling errors} or shortcuts due to applying the circuit laws under invalid assumptions, e.g., use of Kron's reduction in sparsely grounded networks.% Most notably: \textit{a)} use of Kron's reduction of neutral in sparsely grounded networks and \textit{b)} transposition in unbalanced systems. %\begin{enumerate}
       %\item using Kron's reduction (KR) of neutral in networks where they aren't pervasively grounded;
       %\item assuming transposition when there is none.
   %\end{enumerate}
    \item \emph{network data errors}, e.g., wrong topology information;
    \item (largely) inevitable \emph{measurement errors}: `noisy' or `bad' data due to sensor tolerances or malfunction; and
    %\item measurement synchronization error, e.g., assuming everything is measured instantaneously;
    \item \emph{measurement inadequacy} due to 1) \textit{semantic} mismatch, e.g., having averaged instead of instantaneous rms values, %or rms magnitudes instead of base frequency magnitudes; 
    2) \textit{granularity} mismatch, e.g., aggregated three-phase  measurements instead of per-phase; 3) \emph{label} mismatch, e.g., wrong location or phase meta-data, 4) measurements that are not adequately synchronized. 
    % \item nonlinear  equations can lead to many solutions, only one representing reality; this problem can sometimes be avoided with good initial conditions.
\end{itemize}

\noindent Inherently, \emph{any} error adds to the SE residuals, and distinguishing model error, network data error and measurement error is nontrivial. Potential solutions include:
\begin{itemize}[leftmargin=*]
    \item flexible models with variables to represent alternative states, e.g., neutral conductor breakage;
    \item exploiting time-independence.%, so that explicit data isn't required, e.g. impedance.
\end{itemize}
%Other modeling errors include, volt-var control phase-to-ground instead of phase-to-neutral or phase-to-phase;  when real data is unknown, parametrizing the model with ``most likely" values without statistical considerations.
% When real data is unknown, parametrizing the model with \emph{most likely} values without statistical considerations also leads to errors.
% Furthermore, power-voltage and current-voltage formulations have both been used in the context of PF and SE. In rectangular coordinates, they both result in systems of linear and quadratic nonconvex equations, which are only tractable to solve using non-global search algorithms.
% The power-voltage formulation has spurious 0\,V solutions due to lifting Kirchhoff's current law to power variables \cite{Geth2023}. 
% Excluding such solutions from the feasible set is notably challenging for  four-wire DN models.

% Network model data inaccuracies include:
% \begin{enumerate}
%     \item assuming pervasive grounding of the neutral and performing Kron's reduction for the impedances;
%     \item assuming transposition where there is none;
%     \item modeling volt-var control phase-to-ground instead of phase-to-neutral or phase-to-phase;
%     \item when real data is unknown, parameterizing the model with ``most likely" values without statistical considerations.
% \end{enumerate}

The cost to investigate data  issues by human inspection may be prohibitive, and data-driven methods help decrease that effort. 
However, if only power measurements are available, it is unrealistic to reliably discover network parameters\footnote{Phase identification based on power-only data exists~\cite{AryaMIP}, but is known to present poor accuracy compared to the state of the art.}: voltage is necessary for validation. 
In the presence of network errors, the increased residuals may not reconcile, and some SE methods may diverge~\cite{abur2004book} (chap. 8). %: all error sources above contribute to residuals (as discussed in section~\ref{sec:dqi_demand_pi}), which might not reconcile. 
Nevertheless, DSSE, contrary to PF, can be leveraged to calibrate network models.

\section{State and parameter estimation}\label{sec_state_param_estim}

\begin{figure}[tbp]
  \centering
    \includegraphics[width=1\columnwidth]{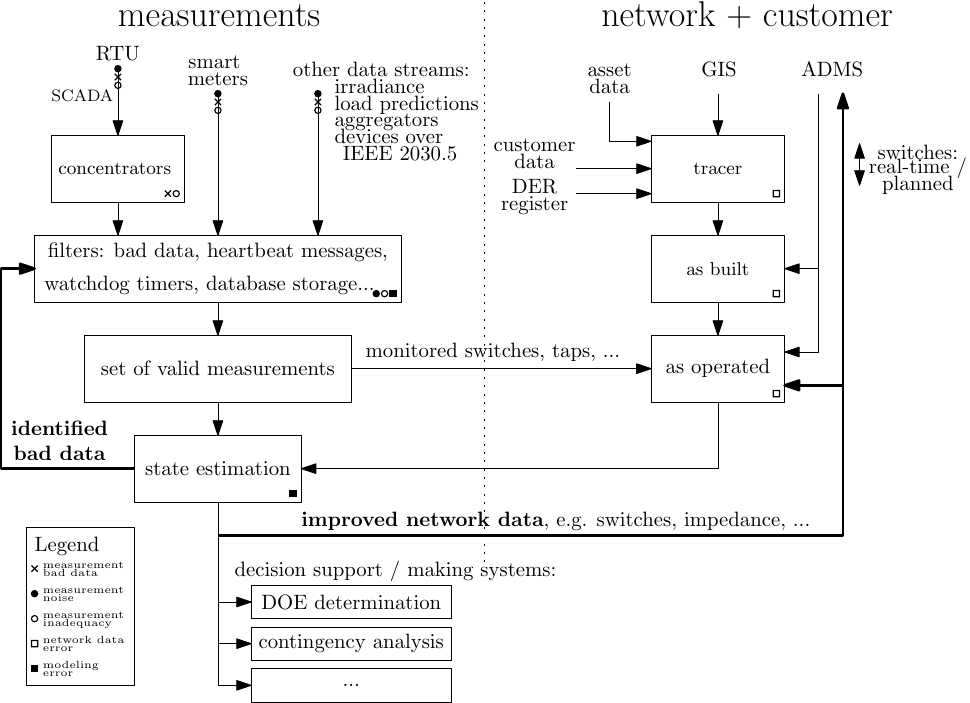}
  \caption{State estimation as part of overall architecture. Novel information flows related to parameter identification indicated bold. }
  \label{fig:se_data_flows}
\end{figure}

\noindent In the presence of pervasive network errors, a first, thorough offline PE is arguably necessary, using historical measurements to derive a DN model as close to the present-state reality as possible. However, this is not sufficient: as power systems change over time, DN data sets should be self-correcting, and such corrections should occur as close to real-time as possible.  
% We note that increases in network utilization can improve impedance estimates, thanks to increasing voltage drops~\cite{Vanin23}. Performing joint state and impedance estimation on \emph{moving time windows} allows to track the quality of the impedance models over time. 
Fig.~\ref{fig:se_data_flows} summarizes our vision for the data flows between measurement and DN data sources towards the SE, indicating the types of data issues that may appear throughout. %from the sources and processing steps throughout. 

From here on, we refer to the first offline PE as `historical' (\apriori), and to the subsequent `on the run' corrections as `near-real-time' (\nearrealtime). Different specifications emerge in the two contexts; e.g., time-independence of properties might be exploitable in  \apriori \  but not in \nearrealtime. 
Table~\ref{tab:categories} provides a summary of parameters to be estimated, including less common ones. 
% Well-known dynamic parameters and typical rates of change are:
% %States/parameters that have more frequently been the subject of study, and their expected rates of change, are the following:
% \begin{itemize}
%     \item Transformer taps: a couple of times per day (online) \footnote{If OLTCs are available, which might be rare in some DNs}, else seasonal, 1x or 2x per year (offline).
%     \item Line codes: consistent for many years.
%     \item Line/cable  impedances: a function of temperature, however variations during certain periods of the year may be very limited. 
%     % Additionally, slow degradations can be observed over the years.
%     \item Transformer impedances: as above, but slower dynamics. %function of temperature, but much slower dynamics than lines/cables. Additionally, slow degradations can be observed over the years.
%     \item Phase connectivity: some power systems perform manual phase switching by congestion, else it is unchanged for years. Automatic switches are still at a research stage.
%     \item Switch states: as discussed in the previous sections.
%     \end{itemize}

%The parameter estimation terminology used in the literature is often inconsistent. Thus, we dedicate the next subsection to shed some clarity, at least in the context of this paper.

\subsection{Notes on Topology Identification}

% The terminology used in the literature is often inconsistent. This section discusses the choices adopted in this paper.

% \subsubsection{Estimation vs identification}

% Semantically, \textit{estimation} suggests an approximate assessment of a certain quantity, and appears better suited when addressing a continuous parameter. \textit{Identification}, on the other hand, suggests a more precise labelling effort, which we believe to be more suited in the context of discrete parameters. We observe that in the literature the terms are used interchangeably,  e.g. in~\cite{Kurup2023, Bariya, Wang2020}. 

% Here we use ``parameter estimation" as generic overarching term to indicate the overall concept. For individual parameters, we attain to the continuous/discrete considerations above. 

% \subsubsection{Topology identification}

% \subsubsection{Topology as a Variable}

\noindent Deka et al.~\cite{DekaTutorial} classify topology learning into topology `detection' and `identification', which is an important distinction as  different methods are used to solve the two. 
The former assesses which lines are energized and which not, assuming their impedances and locations are known, and can be approached with model-based methods.
The latter %is the focus of ~\cite{DekaTutorial}, and 
consists of the green-field reconstruction of the whole connectivity model, (possibly) including line impedances, and is a high-dimensional combinatorial problem. Thus, it is usually addressed with techniques that  abstract from the power flow physics, e.g., graph learning~\cite{Bariya, Pengwah2021}.
%We note that past literature uses the wording ``topology identification" to indicate both concepts, e.g.,~\cite{Karimi2021} for the former and~\cite{Pengwah2021} for the latter. 
% As both are used in the literature, the terminology choice in~\cite{DekaTutorial} is somewhat arbitrary, but 
\emph{
This paper primarily focuses on model (power flow equations) based solutions. 
Thus, we focus on TI in its detection declination.
} 
Interested readers can refer to~\cite{DekaTutorial} for a thorough survey of TI methods. 
% The \emph{nodal admittance matrix} offers a solid basis to look at uniqueness properties \cite{9858017}.
% However, while theorically useful, it is insufficient when looking at the problem of cleaning up existing network data sets \emph{as they are stored in the utility databases}, as nodal admittance matrices have to be rebuilt when the network changes.

%Green-field topology reconstruction is a highly dimensional combinatorial problem, and is not that interesting from the perspective of incremental cleaning of network models. In the literature, green-field reconstruction is usually addressed with techniques that can abstract from the power flow physics, e.g., using graphical learning~\cite{Bariya} \textcolor{blue}{add others from existing list!}.

\subsection{State and Parameter Estimation Literature}\label{sec:state_parameter_lit}

\noindent Table~\ref{tab:parameter_literature} provides a selection of PE references. This is not exhaustive: the  objective of this paper is to discuss the main  ideas behind the methods, not to provide a complete review. 

% W.r.t. Table~\ref{tab:categories} we leave out ``switch state estimation" as it is a sub-category of topology detection.
% Line type and temperature estimation are sub-categories of impedance estimation. The former estimates impedance values and maps them to the closest matching line code from a list. The latter assumes a known nominal impedance and estimates temperature-induced deviations.  

%Carson's equations are commonly solved with geometry (mutual distances between conductors), material-specific impedance at operational temperature and frequency, cross sections and geometric mean radius as inputs, to obtain the capacitance and the series impedance matrices as outputs~\cite{Kersting2011}. 
% The space of the feasible impedance matrices can therefore be reduced when we use information on the inputs to Carson's equations, e.g. known conductor material, lower and upper bounds on the cross-sections and/or distances between the conductors per the installation requirements, spatial configurations such as in-a-horizontal-plane or in-a-triangle. Note that those same inputs are also used in the derivation of the shunt admittance.

\emph{State and parameter estimation can be performed sequentially, jointly, or disjointly}. 

The `\textbf{sequential}' approach was the first to be proposed in the context of integrating PE as part of static SE workflows (although still used in recent literature~\cite{Khalili2022}), and is based on the statistical analysis of SE residuals. First, SE is performed (on one snapshot). Then, the residuals are analyzed in search of \emph{anomalies} that can be attributed to unreported changes or other localized errors, e.g., a wrong line impedance. This concept has been applied %to bad data identification and to the estimation of various parameters (see Chap. 7 and 8 of ~\cite{abur2004book}), and 
particularly in a \nearrealtime \ context, as $1)$ SE is originally designed as a real-time operational tool, and measurement history is not stored, and $2)$ in \nearrealtime \ only few errors are expected (e.g.,  few switch state changes): a diffused presence of errors would pervasively increase residuals, compromising their statistical analysis. In equality-constrained SE, Lagrangian multipliers can replace residuals. %While relatively mature, sequential approaches are still adopted in recent literature~\cite{Khalili2022}.

The joint presence of bad data and network errors is also likely to compromise residual analysis. 
The NYSERDA report~\cite{nyserda} identifies the distinction between bad data and inaccurate network models as a key challenge to real-life DSSE. Measurement scarcity is an additional challenge~\cite{nyserda,krause2015under}. 
Multi-period problems are known to balance lower measurement redundancy levels~\cite{abur2004book} (Section 7.7), and thus appear more promising for DSSE.
Comparing the residuals of separate successive SE calculations, without setting up a real multi-period problem has also been investigated~\cite{Donmez2022, Lin2017}. 
%These do not assume time-independence of parameters, on the contrary, they aim to identify which parameters changed over subsequent measurement \emph{scans}. 
Multi-period implementations are possible with both joint and disjoint approaches, and produce \emph{automated} results. %The incorporation of multiple time periods in a single problem can be pursued with joint and disjoint state and parameter estimation, exploiting the time invariance of the unknown parameters, and producing more ``automated" results. This is possible with joint and disjoint approaches. 

`\textbf{Joint}' state and parameter estimation has existed for several decades for continuous parameters, and is also known as `state (vector) augmentation'~\cite{abur2004book}. Over the years, multiple parameters has been incorporated: multi-conductor impedances~\cite{Vanin23,Mojumdar2021}, zero-impedance branches~\cite{Monticelli1991}, switch~\cite{Karimi2021,Soltani2023} or breaker~\cite{Kekatos2012} states,  tap settings, modeled either continuous~\cite{Korres2004, Therrien2013} (with rounding), or  discrete~\cite{Nanchian2017}. 
% Contrary to other parameters, the impedances in~\cite{Mojumdar2021, Vanin23} do not change over time so the estimation can be done offline. %no point in dong it online overburdening SE. They do it on multi-snapshots, and embedded in state estimation ( as equality constraints).

\begin{table*}[b]
\vspace*{-\baselineskip}
\caption{Parameter estimation literature. \notaddr: not addressed, \textbf{as far as the authors know}. \textcolor{red}{Red} indicates multi-period problems.}
\centering
\begin{tabular}{ll|cc|cc|cc}
\hline
Category & Parameter & \multicolumn{2}{c|}{Joint} & \multicolumn{2}{c|}{Disjoint} & \multicolumn{2}{c}{Sequential} \\
   &      & \multicolumn{1}{c|}{\apriori}    & \nearrealtime    & \multicolumn{1}{c|}{\apriori}  & \nearrealtime  & \multicolumn{1}{c|}{\apriori}  & \nearrealtime  \\
         \hline
Topology & Topology detection & \notaddr & \cite{Karimi2021, Fernandes2023} & \textcolor{red}{\cite{He2021}} & \cite{Gotti2021} & \notaddr & \cite{abur2004book, Lourenco2015} \\
& Topology estimation (green-field) & \notaddr & \notaddr & \textcolor{red}{\cite{DekaTutorial}} & \textcolor{red}{\cite{Babakmehr2016}} & \notaddr & \notaddr \\
& Topology estimation (incremental) & \notaddr & \notaddr & \notaddr & \notaddr & \notaddr & \notaddr \\
& Phase connectivity identification & \textcolor{red}{\cite{Vanin_phaseid}}  & \textcolor{red}{\cite{Vanin_phaseid}} & \textcolor{red}{\cite{Blakely2019, HOOGSTEYN2022}} & \notaddr & \notaddr & \notaddr \\
& Meter-to-transformer assignment & \notaddr & \notaddr & \textcolor{red}{\cite{Hu2023, Saleem2023}} & \notaddr & \notaddr & \notaddr \\
\hline
Line & Impedance estimation & \textcolor{red}{\cite{Vanin23}} & \notaddr & \textcolor{red}{\cite{ClaeysCIRED2021, Zhang2021}} & \notaddr & \notaddr & \cite{abur2004book, Lin2017}\\
& Line type estimation & \notaddr & \notaddr & \textcolor{red}{\cite{Cunha2020}} & \notaddr & \notaddr & \notaddr \\
& Line properties estimation & \notaddr & \notaddr & \notaddr & \cite{Moutis2023} & \notaddr & \notaddr
 \\
% & Carson's equations estimation & \notaddr & \notaddr & \notaddr & \notaddr & \notaddr & \notaddr \\
 \hline
Transformer & Impedance estimation & \notaddr & \cite{Dobakhshari2021} & \notaddr & \cite{ZhangPESGM2015, Aravind2021} & \notaddr & \notaddr \\ 
& Transformer tap position & \textcolor{red}{\cite{Mojumdar2021}} & \cite{Nanchian2017, Korres2004, Therrien2013} & \notaddr & \notaddr & \notaddr & \cite{Pires2014} \\ 
& Transformer type & \notaddr  &  \notaddr  &  \notaddr  & \notaddr  & \notaddr &  \notaddr  \\
 \hline
Regulator/OLTC & Gang vs per-phase operation & \notaddr & \notaddr & \textcolor{red}{\cite{yusuf2022}} & \notaddr  & \notaddr & \notaddr \\
& Other (impedance, control...) & \notaddr  &  \notaddr  &  \notaddr  &  \notaddr  & \notaddr &  \notaddr \\ 
 % &  Controller mode estimation &  \notaddr  & \notaddr & \notaddr & \notaddr  & \notaddr & \notaddr 
 % \\
 % & Bidirectional capable? & \notaddr  &  \notaddr  &  \notaddr  &  \notaddr  & \notaddr &  \notaddr  \\
  \hline
PV system & PV rating estimation & \notaddr & \notaddr & \textcolor{red}{\cite{Mason2020}} & \notaddr & \notaddr & \notaddr \\
 & PV control mode identification & \notaddr & \notaddr & \notaddr & \notaddr & \notaddr & \notaddr 
 \\
 & PV control setting estimation & \notaddr & \notaddr & \textcolor{red}{\cite{Talkington2022}} & \notaddr & \notaddr & \notaddr \\
 \hline
\end{tabular}
\label{tab:parameter_literature}
\vspace*{-1.5\baselineskip}
\end{table*}

Augmenting the state to incorporate discrete parameters would lead to mixed-integer non-convex programming problems, which are prone to scalability issues. 
Relaxations and approximations of non-convex constraints have led to tractable yet effective mixed-integer problems. 
These are relatively common in the context of joint SE and topology detection, e.g.,~\cite{Karimi2021}, but are nearly non-existent for phase identification: deriving the connectivity of \textit{all} users results in a large number of binary variables, and to our knowledge, the only mixed-integer SE augmentation for phase connections is~\cite{Vanin_phaseid}. 
%On the other hand, in phase identification, the connectivity of \textit{all} users is derived. This results in a large number of binary variables, and to our knowledge, the only mixed-integer SE augmentation for phase connections is~\cite{Vanin_phaseid}. 

% Note that augmented DSSE that preserves exact non-convex equations, e.g.,~\cite{Vanin23}, returns the most-likely parameters and state values, in principle without approximation. When linearizations or relaxations are applied, modeling errors are introduced, and, in most practical settings, a re-run of exact DSSE after fixing  the derived parameters would be desirable. %\footnote{Nevertheless, augmented DSSE has advantages from a pure parameter estimation standpoint, e.g., in terms of amount of measurements needed. See discussion in~\cite{Vanin_phaseid}}.

While all modeling errors usually increase the residuals, inaccurate impedances generally cause much lower increases than wrong discrete parameter values~\cite{abur2004book, Vanin_phaseid}. This can be exploited, e.g., first to estimate discrete variables, fix them, and then subsequently perform impedance estimation~\cite{Vanin_phaseid}.

% Finally, while these are not network properties, DSSE can be augmented with power variables for unmonitored users to perform load allocation. 
% This could be particularly useful for three-phase users that present only total power measurements, that need to be subdivided across the three phases.

`\textbf{Disjoint}': PE is performed without using SE concepts, usually in an \apriori \ manner. Phase identification is almost exclusively performed disjointly, most often via voltage measurements clustering~\cite{Blakely2019}, although power measurements can be used~\cite{HOOGSTEYN2022}, and the same holds for green-field TI~\cite{DekaTutorial}. Several researchers estimate cable impedances with linear regression~\cite{Lave2019, Peppanen2016, Cunha2020}. Yusuf et al.~\cite{yusuf2022} present data-driven methods  
% \footnote{including statistical analysis of voltage measurements and similarity maximization.} 
to identify the tap position of voltage regulators and establish whether the tap changers are phase- or gang-operated. %Green-field TI also belongs to this category and is generally performed on historical data~\cite{DekaTutorial}. 
Gotti et al.~\cite{Gotti2021} and He et al.~\cite{He2021} use DSSE to validate topologies estimated in a disjoint manner. % add that they are in a sense sequential but actually disjoint. NN is ofc trained offline. In a sense, all other topology reconstruction would also solve the detection problem from the apriori perspective.

%\cite{He2021} does use SE to validate between all possible topologies a priori

Note that topology detection is usually limited to detecting switch states when it is done jointly with or sequentially after DSSE. % A more generic selection of energized lines given a subset of possible topologies has been investigated in~\cite{Pengwah2021} in a disjoint manner.
TI methods (or adaptations) could perform `incremental TI', i.e., sparse connectivity errors that do not correspond to switch locations. These data quality issues exist in real network data~\cite{GethCIRED2023} but are not addressed in the literature.
Topology, line parameters, and tap settings appear well investigated, whereas limited work addresses other transformer/regulator parameters. 
In particular, transformer types and groups are assumed to be known, while in our industrial experience this is not the case. 
Transformer models are defined by winding resistances, leakage inductance, shunt properties, etc. 
Data-driven estimation methods exist~\cite{Mossad2014}, which do not require  disconnection. 
% There are similarities between these methods and DSSE, e.g., minimization of residuals and use of least square regression. 
However, transformer parameter estimation is performed with `local' measurements only, e.g., voltages and currents at the transformer terminals, and not as part of a `full' augmented DSSE, except in~\cite{Dobakhshari2021}. However, \cite{Dobakhshari2021} only estimates turns ratio and series impedance. 
A few references are gathered in Table~\ref{tab:parameter_literature} as `Transformer - impedance estimation'. Notably, most papers address single-phase transformers only.

Finally, most - if not all - existing PE methods implicitly assume that the users/meters are assigned to the correct transformer (beforehand). If this is not the case, the successive estimation process is compromised. Researchers have only recently addressed the data-driven association of meters to the correct transformers. While the literature is limited, various names are used, e.g., `transformer-customer relationship identification'~\cite{Hu2023}, or `meter-transformer mapping'~\cite{Saleem2023} (M2T). 
% In North-American-like secondary circuits, where transformers serve one or few single-phase users, M2T resembles phase identification. To the authors' knowledge, M2T has not been applied to DNs with multi-phase secondary. This is a very relevant area for future research.

%Finally, both DSSE and purely data-driven methods are developed in~\cite{osti_1888157} (chap 8) to track voltage regulator and capacitor states.

\subsection{Behind-The-Meter Resource Estimation}

% To simulate a network, PV inverter controllers, and tap changer controllers are necessary inputs to obtain realistic results, whereas in SE the set points of these controllers can be resolved without \emph{requiring} controller models.
% Identifying controller models is quickly becoming a major challenge for scenario analysis.
% For instance in LV networks, PV inverters with a variety of configurations (constant power factor vs Volt-var/Watt) and parameterizations are deployed simultaneously. 
% Weaker networks behave quite differently when dominated by different kinds of inverters.

\noindent \emph{Parameter estimation}  typically refers to network asset properties, however, parameters of behind-the-meter (BTM) resources are often unknown 
% \footnote{Especially as they are not owned or accessible to utilities.}, 
and assuming certain properties may significantly impact grid analytics and monitoring~\cite{IEEETFDER}. The identification of BTM parameters is a novel research topic and has not been explored in the context of augmented SE, even though for some parameters this might be possible. 

In contexts where smart meters (SMs) are not bidirectional, disaggregation of demand and PV generation is essential. Existing approaches are reviewed in~\cite{IEEETFDER}. Like for DN PE, these can be model-based, fully data-driven, or a combination.% in-between. %in between (physics-informed data-driven). 

Where SMs are bidirectional, grid injections may be directly measured (although self-consumption may still be unaccounted for), but inverter control settings and energy management behavior are still usually unknown. %, especially at the grid edge. 

Most methods for PV require or are facilitated by the use of exogenous information, such as irradiance~\cite{IEEETFDER}, although methods based only on electrical measurements exist~\cite{Mason2020, Talkington2022}. Mason et al.~\cite{Mason2020} present a supervised learning model for size, azimuth and tilt estimation of PV installations based on net power data. Talkington et al.~\cite{Talkington2022} disaggregate demand and generation and derive reactive power control settings (constant PF and Volt-var) of PV inverters, using historical power and voltage data. Voltage-only-based methods are also possible, if a good-quality secondary circuit model is available~\cite{osti_1888157}. 

There is significantly more literature on PV disaggregation than control mode estimation, e.g.,~\cite{Bu2023}. Similarly, the non-intrusive load modelling literature is also abundant~\cite{IEEETFDER}, and techniques similar to~\cite{Talkington2022} could be extended to other inverter-based resources, post disaggregation of the resource's load. To our knowledge, there is no literature on the topic. 

% In contexts where users are not required to notify IBR installations to their utility, the data-driven localization of such resources is relevant. See~\cite{osti_1888157} for the case of PV.

Additionally, the identification of PV control modes appears unaddressed, e.g.,~\cite{Talkington2022} assumes to know a-priori which inverters perform constant power flow control and which Volt-var.  
% There appears to be room to build on this work, and also explore other control strategies, e.g., Volt-Watt.
Similarly, identifying regulator control modes and bidirectional capabilities is industrially relevant but unexplored. 

Finally, we note that the identification of load models (constant power, vs ZIP, etc.) is also beneficial, and has received academic attention~\cite{IEEETFDER}. 

% \subsection{Carson}

 % In the table, line type and temperature estimation are sub-categories of impedance estimation (IE). The former estimates impedance values and maps them to the closest matching line code from a list. The latter assumes a known nominal impedance and estimates temperature-induced deviations.
 % Carson's equations could also be used to calculate impedances as functions of conductor properties: cross-section, material, etc.~\cite{Kersting2011}. Their use could reduce the search space of the IE problem, but has not been investigated yet.
 
% In the context of SE, load models considerations are usually ignored, but their incorporation would result in, albeit perhaps slightly, different estimations (see discussion in section~\ref{sec:beyond_canonical_state_estimation}).  

\section{Necessary Algorithmic Advancements}\label{sec_algorithms}

\noindent We now discuss algorithmic bottlenecks of  DSSE generalizations. 
We posit a need to move beyond Gauss-Newton (GN) and WLS when the following are desired:
\begin{itemize}[leftmargin=*]
    \item non-Gaussian state estimation (both in a robust sense and in a best estimator sense);
    \item conditional estimation;
    \item joint state and parameter estimation;
\end{itemize}
We argue that all are crucial to make DSSE practical.

\subsection{Gauss-Newton and Other Algorithms}
\noindent \textit{Newton-Type WLS Algorithms:}
The most common solution technique for SE is the GN algorithm \cite{monticelli2000electric, Wang2018} (\S V.A), a Quasi-Newton method that is applicable to least squares minimization problems. GN approximates the Hessian from the Newton-Rhapson (NR) method in its iterations, making it less robust than NR. Furthermore, it cannot naturally incorporate variable bounds or inequalities. Nonetheless, by avoiding the  calculation of the Hessian, GN is generally assumed to require less computation per iteration than NR.
% It minimizes the least squares objective by taking iterates, approximating the Hessian from the full-Newton method.
%Because GN approximates the Hessian value, it is less robust than the full Newton approach and cannot be naturally applied to constrained optimization problems with bounds and inequality constraints.
% Nonetheless, because it does not require full inversion and calculation of the Hessian, it can be computationally faster than full-Newton.

The incorporation of very diverse weights for different measurement types in DSSE, particularly pseudo-measurements (very low weights) and  zero-injection buses (very high) results in  ill-conditioning. Improvements on the basic GN-based implementation avoid using high weights for virtual measurements~\cite{dzafic2017}, modelling them as equality constraints and solving the Lagrangian~\cite{monticelli2000electric}.  Further extensions of these algorithms allow for SE in underobserved systems, where the gain matrix is rank-deficient and thus non-invertible~\cite{krause2015under}.

% For SE, GN necessitates using concepts such as `state vector', `zero injection nodes' and `pseudo measurements'.
% Initially, the GN-based SE approaches included zero-injection nodes through the use of high weights; however, it resulted in ill-conditioning.
% More advanced versions of the GN-based implementation avoided using high weights for zero-injection buses \cite{dzafic2017} in equality-constrained WLS estimators \cite{monticelli2000electric}. 
% Further extensions allow for state estimation in underobserved contexts \cite{krause2015under}. 

NR has also been explored for constrained optimization-based SE, but much less frequently. 
 NR can include any number of equality and inequalities, and (with heuristics) has better convergence for large-practical networks, which helps solve difficult instances. In practice, NR convergence significantly improves if \emph{i)} the Hessian at each iteration is positive-definite, to ensure descent direction and \emph{ii)} LS techniques are used to choose an appropriate step size to achieve sufficient objective improvement in each step. Due to its approximations, GN does not need routines to deal with non-PSD Hessians, but, as a consequence, LS may not be as powerful.

\textit{Trust-region  algorithms:}
The Levenberg-Marquardt (LM) method, a kind of trust region method that applies to least-square minimization problems, improves on GN. LM interpolates between the GN and the method of gradient descent and is more robust than GN: it can generally find a solution even when the initial condition is far from the minimum. Moreover, it is more stable than NR as it ensures step-sizes within a trust region but generally converges slower than GN and NR.

Very little literature exists on trust region-based SE.
Pajic and Clements \cite{Pajic2005} explore them for TSSE, together with GN with line search (LS), motivated specifically by improving reliability in the presence of bad data such as topology errors.
They observe that the trust region approach is the most reliable in difficult conditions, e.g., the presence of cascade failures.

% Discuss Quasi-Newton and Levenberg Marquardt? and gradient based algorithm\cite{Zhou2020Gradient}?

% Advancements on the OR front are faster than those in the power systems community.

\textit{Nonlinear programming:}
Recently, NLP approaches have become  popular \cite{Melo2019, Melo2023, Singh94, dzafic2014real, li_circuit_se, li2021convex_gse}. 
%Algorithms for equality and inequality constrained optimization can be used to solve SE problems. 
These can solve SE with both equality and inequality constraints. 
Inequalities can be used to model non-monitored users~\cite{Melo2019}. 
%Melo et al.~\cite{Melo2019, Melo2023}, for instance, model non-monitored users as inequality constraints. 
Algebraic modeling toolboxes with automatic differentiation  \textit{1)} allow a clear separation between the model ( equations and  objective) and the solver (e.g. \textsc{Ipopt}), \textit{2)} can support WLAV objectives naturally \cite{li2021wlav, Singh94}, as the nondifferentiability of the objective function is avoided through reformulation using inequalities and the epigraph transform, and \textit{3)} facilitate alternative formulations of the circuit physics, e.g. through convex relaxations~\cite{Zhang2018} or linear approximations.
% Approximations of the physics, e.g. through convex relaxations (unbalanced branch flow model SDP relaxation) \cite{Zhang2018} or linear approximations (e.g. LinDist3Flow) have also been explored.
There is much less research on alternative formulations in the context of SE than OPF.  Presumably, this is because any approximation implies modeling errors that add up to the residuals, making the decision support outcome measurably worse despite potential speedups. For PE, approximations may be acceptable, but when approximated SE and PE are performed jointly, a final nonconvex SE with the fixed parameters is recommended. 
Taheri et al. discuss using  SE to recover feasibility of relaxed/approximated OPF models~\cite{Taheri2023}.

SE can also be approached with matrix completion, which consists of estimating missing values in low-rank matrices. 
Such methods are particularly well-suited for contexts with very limited observability \cite{Donti2020}, but most implementations rely on semidefinite programming (SDP) models \cite{Hu2022, Liu2020, Donti2020, Rout2022}, which -- despite their convexity -- struggle to scale in practice. Specifically, the network feasibility in the SDP approach is inherently dependent on the solution being rank-1; which is rarely the case for larger networks. If the solution is not rank-1, a feasible AC solution cannot be recovered. % from the lifted SDP solution.
Therefore, Liu et al. develop a more practical approach as a quadratic nonconvex optimization problem \cite{Liu2019}. 

% Tensor completion \cite{Liu2022}.

\textit{Other approaches:}
Cao et al. \cite{Cao2023} give an overview of learning-based DSSE methods and propose one that can detect topological changes by including some physics. 
However, the rationale behind adopting a learning-based approach is that system parameters are unknown, which we argue to be an issue that should and can be solved instead.
% Wang et al. \cite{Wang2018} propose a ``feasible point pursuit'' approach. 
% Alternatively, Pau and Tamim \cite{Pau2023} use a variant on the backward-forward sweep (power flow) method, that is able to resolve measurements including voltage. 

\textit{Discussion:}
% Advancements on the OR front are faster than those in the power systems community.
% Robust estimators generally aim to efficiently identify bad data when the error distribution is unknown. When the latter is known, more generic MLE frameworks can be set-up, to always use the optimal estimator even in non-Gaussian conditions~\cite{Cheng2021, VaninMLE2023}. 
% We believe that such approaches are rarely considered because: \textit{1)} it is not solvable with Gauss-Newton, and \textit{2)} the DSSE literature somewhat ``drifted away" from the mathematical foundations, towards the narrow use of WLS methods.
We posit that for large problems, especially those with multi-period constraints and many zero-injection nodes, NR approaches within constrained optimization SE frameworks are the more robust, i.e., most likely to converge.

Approaches that support inequality constraints allow for \emph{conditional} SE: e.g., to establish a quality of the fit while requiring variables to be in a certain range. 
Furthermore, in GN, the objective must be a WLS function or variants such as a Schweppe-Huber. WLAV, however, is nondifferentiable, so it requires a different algorithm.
% , which is interesting because it is linear and robust to bad data is not possible as it is nondifferentiable (without reformulation).

% \subsection{Variable Spaces for Describing the Physics}

%In the context of PF and SE, both the (rectangular) power-voltage and current-voltage formulations result in systems of linear and quadratic nonconvex equations.
% However, if we generalize the problem statement to  variable impedances, IVR remains quadratic, whereas SVR becomes cubic \cite{Vanin23}\footnote{Obviously, we can re-cast this as a system of quadratic equations by introducing auxiliary variables to represent certain bilinear terms.}. 

\subsection{Maximum Likelihood Estimation}
% Mathematically, WLS SE is a maximum-likelihood estimator (MLE) only if all measurements present white Gaussian noise. 
% More generally, SE is a MLE problem, where we try to establish realized values for random variables subject to a probability distribution. 
% Different distributions have different maximum likelihood models \cite{VaninMLE2023}, for instance:
\noindent SE is a maximum likelihood estimation (MLE) problem that establishes realized values for random variables subject to a probability distribution. WLS minimization is the `correct' MLE only if all measurements present white Gaussian noise. WLS corresponds to the $\ell_2$ norm, with strong theory and closed-form results.  Different distributions have different maximum likelihood models \cite{VaninMLE2023}; for instance, the Laplace distribution maps to the $\ell_1$ norm, i.e., WLAV minimization.
% \begin{itemize}
%     \item Gaussian uncertainty maps to the $\ell_2$ norm a.k.a. WLS, where there is strong theory and closed-form results;
%     \item the Laplace distribution maps to the $\ell_1$ norm a.k.a. WLAV.
% \end{itemize}
As a regularizer, $\ell_1$ is understood to induce sparsity much more than $\ell_2$, which makes the WLAV more robust. WLS is sensitive to outliers w.r.t. to alternatives, so it is sometimes deliberately replaced with  objectives that are more robust:
\begin{itemize}[leftmargin=*]
    % \item WLS SE is sensitive to outliers relative to alternatives, so it is sometimes deliberately not chosen;
    \item $\ell_0$ `norm', a.k.a. measurement discarding, may identify sources of bad data but is combinatorial;
    \item the Schweppe-Huber function avoids nondifferentiability of WLAV, is less sensitive to outliers but  typically nonconvex;
    \item the matrix nuclear norm as a proxy for rank in the context of matrix completion \cite{Donti2020} is convex;
    % \item likelihood models, including non-Gaussian, 
    \item mixing the above, with or without penalization.
\end{itemize}

% Robust estimators are one way to handle bad measurement data, which consist of extremely unlikely (high sigma) events that defy the used probabilistic error characterization, and thus result in distorted estimates.  However, bad data rejection is still a major challenge for DSSE.

% Singh et al discuss how to pick the right fitting criteria \cite{singh2009choice}

% \begin{itemize}
%     \item ivr more reliable
%     \item ivr remains quadratic for impedance estimation, svr does not
%     \item ivr is necessary for 4-wire
%     \item problems with spurious solutions
% \end{itemize}

% \subsection{Error Rejection Strategies}
%  In this context, we note the lack of generally available statistical models to pick up on outliers. 

\section{Real-World Experiences with DSSE}  \label{sec_deployment}
\subsection{Lessons Learnt From Collaborative R\&D Projects}
\noindent Bad data can only partially be managed through heartbeat protocols, watchdog timers and error correction codes.
In addition to large random noise events, data \emph{adequacy issues} and flawed sensor \emph{installations} are common in practice:
% The following are particularly detrimental:
\begin{itemize}[leftmargin=*]
    \item unknown or inconsistent measurement setup, e.g., phase-to-neutral vs phase-to-ground voltages,
    \item miswired voltage/current transformers (i.e. 180$^{\circ}$ offset);
    \item wrongly accounted-for voltage/current transformer ratios; % not being taken into account;
    \item phase rotation order wrong;
    \item mislabeled units (e.g. kW vs W, V vs kV);
    \item measurement values reported are min/max over the interval, not the mean or median.
\end{itemize}
Some of these errors are expected to be caught  promptly, e.g., when being used for billing purposes. 
% Nevertheless, issues may persist, even if the frequency of occurence is low, when additional sensors are installed.  
% Different information systems may handle technical data vs billing data; and corrections may happen at any stage throughout the information flow. 
% However, , such issues  lead to  delays. 
During SE deployment, experimentation with  error rejection strategies can be particularly useful as a mechanism for root-cause analysis. 
% These are discussed in~\ref{sec:root-cause-analysis}, where the benefit of performing inequality-constrained DSSE emerges, which calls for the NLP strategies discussed earlier in this section.
Finally, we observe that while they may not result in bad data, the following adequacy issues are detrimental to SE and/or PE:
\begin{itemize}[leftmargin=*]
    \item long measurement intervals, thereby underestimating spiky loads and network losses;
    \item only aggregate three-phase power measurement available, not individual phase or line values;
    \item measurement systems that only communicate new values when a (high) threshold for change is exceeded.
    % \footnote{this compounds with a missing heartbeat protocol}.
\end{itemize}

\subsection{Lessons Learnt At Alliander DNO}
\noindent Alliander, the largest DN operator in The Netherlands, serving over three million customers, started in 2015 with the implementation of SE.
Its MV and LV network consists of over 40\,000 km of cable split up into 22 million segments. 
% Alliander operates both the MV and LV network.

The biggest technical challenges were developing a fast enough DSSE algorithm and preparing a suitable dataset for the entire DN. At the time, no marketed product could handle SE of this scale, so Alliander initiated the development of the \textsc{PowerGridModel}~\cite{Xiang_PowerGridModel_power-grid-model}, which according to Xiang et al.~\cite{Xiang2023}, the fastest open-source SE toolbox currently available.
Alliander follows a SE implementation pathway, with foreseen applications in increasing order of practical difficulty:

\textbf{Improving topological network data}: 
The first result of every practical SE attempt is finding errors in topological network data. In the case of Alliander's dataset of 22 million cables, hundreds of thousands of errors were identified. 
    The faults mainly consisted of short, missing  LV cable sections, creating network areas seemingly disconnected to the main grid. Since many processes within a utility rely on accurate data, this is already a very valuable result.
    Most SE algorithms also give spatial information about mismatches between models and sensor information. This helps identify faulty sensors and unrealistic cable properties.

 \textbf{Optimizing sensor design and placement}: Certain SE implementations can give the uncertainty of their estimated states. This is valuable information for adding sensors to the network, as they can be added to the places with the most uncertainty. Using simulations, it can also be estimated how much extra sensors improve the results of SE.

    \textbf{Real-time load modeling MV}: Once data and sensoring are sufficient, SE is deployed, generating valuable information on network congestion and outages. 
    It is advisable to start the practical implementation of SE on the MV network, as the LV network poses extra challenges. % regarding data quality and computation.

    \textbf{Real-time congestion management}: With clear real-time insight of the energy flows within the electricity network, techniques to control congestion can be explored. 

    \textbf{Fraud detection}: With near-perfect SE, detecting fraud is a possibility. However, given the cost of a false positive, it is probably infeasible to create an SE result that is good enough to use for fraud detection on its own.

The biggest non-technical challenge of implementing SE was keeping the algorithms explainable to non-technical colleagues. 
A metaphor that worked well is, `We are currently regulating traffic without knowing where the traffic is. SE will let us know where the traffic is in real-time.'. Alliander believes SE will be invaluable in finding and mitigating congestion, thus facilitating the energy transition.

\subsection{Lessons Learnt At Energy Queensland}
Two DNs operate in the state of Queensland, Australia:
Energex services the heavily populated southeast corner, and Ergon Energy Network the regional and remote northern parts. 
They connect 2.3 million customers with 212\,000\,km of powerlines and cables ranging from 132\,kV to LV. 
Queensland DNs support almost 745\,000 rooftop and other small-to-medium solar energy systems, with numbers steadily increasing. 
% The proliferation of rooftop solar
This triggered minimum demand situations, significant reverse flows, voltage performance issues and DN congestion. 
In 2019, these challenges prompted the implementation of DSSE to improve DN visibility using the limited telemetry available. 
Key lessons learned are:

 \textbf{Network Model:} The term `network model' means many things to many people. 
    Network models used in other business contexts, e.g., operations and planning, shortcuts for simplifying visual representation, computation or hold only hierarchical associations, resulting in deficiencies in the detail required by DSSE. 
    The Geospatial Information System (GIS), being the only source of LV network data, was deemed the primary network model data source for SE; however, the GIS is filled with connectivity issues, particularly at LV. 
    One example is that customers are represented as being attached to a distribution transformer rather than the actual point of coupling to the LV DN. 
    Improving field data capture and using SMs to infer connectivity improves DN data. 
    MV is the initial focus of DSSE deployment due to the superior quality of available network models w.r.t. LV.

\begin{figure}[tbp]
  \centering
    \includegraphics[width=1\columnwidth]{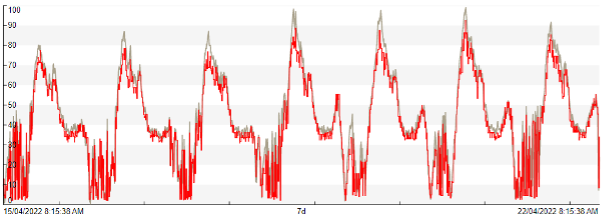}
  \caption{Out-of-sample validation of underobserved DSSE on an MV DN feeder in Queensland. The indicated quantity is current (A), the grey line is the DSSE solution, and the red line is measured. }
  \label{fig:validation}
\end{figure}

    \textbf{Telemetry:} 
    SMs are the responsibility of retailers and are used for billing purposes only, making access to SM \emph{power quality} data difficult for DNs. 
    The industry recognises the value of SM data for DN operations, but regulatory changes improving access take time. 
    The rollout of monitoring on the LV side of MV/LV transformers has proved invaluable. 
    Monitoring coverage of $\approx$30\% has generally proven sufficient for accurate SE results (Fig. \ref{fig:validation}). 
    Standards are being updated to improve telemetry from SCADA devices for SE. 
    Limitations include aggregated three-phase power flows rather than individual phase measurements or measuring current magnitude rather than active and reactive power. 
    Setting up data streaming capability to deliver telemetry to DSSE in NRT is fundamental to success and was achieved through Apache Kafka.

    \textbf{Cross-System Data Mapping and Naming Conventions:} 
    SE requires reconciling inputs from GIS, asset, customer, metering and SCADA data systems. 
    Naming conventions are important to ensure objects (e.g., a transformer), asset characteristics (e.g., rating, tap position) and telemetry can be matched. 
    Alignment is challenging when information is mastered in different source systems.

Early results show promise in delivering DSSE-derived `synthetic' grid visibility data as distinguished from sensors. 
Power flows are not measured on most MV feeders making periods of reverse flow hard to identify. 
Currently, the DN utilities use CT data looking for a characteristic bounce at the zero crossing as an indication of likely reverse flow (Fig \ref{fig:validation}). 
% Using SE and the downstream  transformer monitor data, the feeder power flow can be estimated, as shown in Fig \ref{fig:validation} and \ref{fig:reverse}.

% \begin{figure}[tbp]
%   \centering
%     \includegraphics[width=1\columnwidth]{Picture 1.png}
%   \caption{DSSE-derived MV active and reactive power reveals reverse power flows. Grey line indicates active power (kW), purple is reactive power (kvar).}
%   \label{fig:reverse}
% \end{figure}

Another key application of SE and subsequent network-constrained optimisation is to support dynamic connection of customer DERs \cite{9715663}. 
Networks collaborate with customers to manage rooftop solar, electric vehicles and batteries by flexibly varying export and import limits at the point of coupling, informed by congestions made visible through DSSE.
Educating engineers, field technicians and solution architects within the business on incremental no-regrets actions  to establish a digital and data environment favourable to DSSE allows the opportunities to be realised sooner and more cost effectively.

\section{Conclusions and Future work} \label{sec_conclusions}
% We observe that even within this seemingly limited and ``mature" space, there are significant misconceptions. 
% More specifically, we focus on how conventional mechanistic algorithms for DSSE can be improved to address these needs.

\noindent This paper discussed the main technical bottlenecks and industrial misconceptions that hinder the adoption of state estimation in DNs. 
% We share some of the lessons learned in our industrial experience, which academics can benefit from for more ``domain aware" research. 
We argue that treating state and parameter estimation as a single problem, instead of two distinct ones, is a crucial step towards applicable solutions. 
This calls for the adoption of less conventional problem formulations and algorithms. 
Some examples exist in the literature, and have been reviewed here. 
We highlight specific research gaps where present, but stress that, in general, a change of perspective is required, e.g., including adaptive, automated and thorough network model cleaning pipelines. 
% Augmenting the `state' to include parameter variables is necessary, a.o., to facilitate root-cause analysis. Nevertheless, additional solutions in the latter space should and can be pursued. %major benefits can be harvested from the convergence of state and parameter estimation problems , which calls for the the adoption of less conventional models and algorithms. %Subsequently, we analyze solutions 
Finally, further benefits can be harvested through R\&D on:
% in the following, largely underaddressed areas:
\begin{itemize}[leftmargin=*]
    \item exploration of synergies between data-driven and physics-based solutions for state and parameter estimation;
    \item root-cause analysis methods for network (data) issues;
    \item maximum likelihood and sparsity-inducing estimators combined with improved statistical models to detect meter configuration errors and other anomalies; 
    \item implications of limited observability on decision risks.
    % generalized state estimation across multiple time steps, with state variables for taps, switches, etc.
    % \item contingency analysis through power flow acknowledging uncertainty in the state
    % \item data-driven is not in conflict with physics-based
\end{itemize}

\bibliographystyle{IEEEtran}

\end{document}